# AN ALGEBRAIC METHOD FOR OPTIMIZING THE STATE, CONTROL, AND TERMINAL STATE WEIGHT MATRICES FOR OPTIMAL FEEDBACK CONTROL

**Daegyun Choi,[*] Donghoon Kim,[†] and James D. Turner[‡]**

The necessary conditions for formulating optimal feedback control algorithms have been known for many years. Free parameters exist in the performance index in the form of state and control penalty and terminal state penalty matrices for tuning the performance of the optimally controlled system. The selection process is typically experimental and iterative. To generate the weight matrices for optimal feedback control algorithms, an optimization process is proposed. Typically, the optimization process for the weight matrices requires several numerical integration processes that are computationally expensive; this work overcomes the classical high computational cost by exploiting closed-form solutions for the time-varying Riccati matrix, state trajectories, state transition matrix, and the optimal performance index. Closed-form algebraic equations are used to generate all partial derivative calculations, and no numerical integration is required. The closed-form partial derivatives are used to generate analytic gradients for the optimization steps. The optimization strategy seeks to minimize the terminal state values for the feedback control problem. A numerical example is presented to demonstrate the effectiveness of the proposed optimization algorithm. The resulting computational procedures are expected to be broadly useful for control theory applications in science and engineering.

## INTRODUCTION

In the optimal control theory, the optimal control problem is generally defined to determine the optimal control trajectories that satisfy the physical constraints and minimize a quadratic performance index.[1] To obtain the optimal control trajectories, the weight matrices for the state, the control, and the terminal states are iteratively adjusted to tune the performance of the controlled system dynamics. Although a key factor of the optimal feedback control algorithm is to determine the weight matrices properly, these matrices are determined by trial-and-error in most studies, and this is time-consuming work. Many scholars proposed several approaches to reduce the effort for the selection of the weight matrices. Kukreti et al.[2] considered the linearized spacecraft attitude dynamics for the optimal feedback control problem. They applied the genetic algorithm to determine the weight matrices, and the weighted summation of the final time and the attitude and angular velocity error was considered as a fitness function to be minimized. Moreover, including the aforementioned study, plenty of studies utilized the population-based heuristic algorithms, such as genetic algorithm[3,4,5], particle swarm optimization[6], adaptive particle swarm optimization[7], artificial bee colony optimization[8], adaptive predator-prey optimization[9], and Jaya algorithm[10]. Unlike the weight selection methods using the population-based heuristic, some researchers proposed analytical approaches using Lagrange optimization technique[11] and pole-placement approach[12]. In the study investigated by Yang[13], a pole assignment design of a quaternion-based spacecraft control problem was utilized to determine the weight matrices considering the balance between the performance and fuel consumption.

[*] PhD Student, Department of Aerospace Engineering and Engineering Mechanics, University of Cincinnati, Cincinnati, OH 45221, USA.
[†] Assistant Professor, Department of Aerospace Engineering and Engineering Mechanics, University of Cincinnati, Cincinnati, OH 45221, USA.
[‡‡] Retired Aerospace Engineer, AAS Fellow, 9399 Wade Blvd., Frisco TX 75035, USA.

Mostly, each component of the weight matrices in the optimal feedback control problem is optimized by utilizing the population-based heuristic algorithms. In addition, a number of studies determined only diagonal components of the state and control weight matrices and used the algebraic Riccati equation. On the other hand, this work considers the differential Riccati equation and optimizes all the symmetric elements of the weight matrices to provide more flexibility by utilizing analytic gradients, and the optimization strategy seeks to minimize the terminal state values for the feedback control problem. The optimization task is defined by Taylor expanding the state and control at the final time to force the desired final time boundary conditions for the system response. In the optimization process, differential equations are numerically integrated backward and/or forward in time, but those computations may make this approach inefficient. However, this work significantly reduces the computational burden by replacing differential equations that require computationally intensive numerical integration processes with purely algebraic equations of closed-form solutions. To validate the performance of the proposed approaches, the simulation study is performed by applying the proposed approaches for the second-order linear differential system.

## PROBLEM FORMULATION

The optimal feedback control problem is formulated by seeking $\mathbf{u}(t)$ to minimize[1]

$$J = \frac{1}{2}\mathbf{x}^T(t_\mathrm{f})S_\mathrm{f}\mathbf{x}(t_\mathrm{f}) + \frac{1}{2}\int_{t_0}^{t_\mathrm{f}}(\mathbf{x}^T(t)Q\mathbf{x}(t) + \mathbf{u}^T(t)R\mathbf{u}(t))\,\mathrm{d}t \tag{1}$$

subject to

$$\dot{\mathbf{x}}(t) = A\mathbf{x}(t) + B\mathbf{u}(t);\ \mathbf{x}(t_0) = \mathbf{x}_0; \quad \mathbf{x}(t_\mathrm{f}) = \mathbf{x}_\mathrm{f} = \mathbf{0} \tag{2}$$

where $\mathbf{x}(t) \in \Re^n$ denotes the system state vector, $A \in \Re^{n\times n}$ is the system dynamics matrix, $B \in \Re^{n\times m}$ is the control influence matrix, $\mathbf{u}(t) \in \Re^m$ denotes the system control vector, $Q = Q^T \geq 0 \in \Re^{n\times n}$ denotes the state weight matrix, $R = R^T > 0 \in \Re^{m\times m}$ denotes the control weight matrix, and $S_\mathrm{f} = S_\mathrm{f}^T \geq 0 \in \Re^{n\times n}$ denotes the terminal state weight matrix. The optimal control is given by

$$\mathbf{u}(t) = -R^{-1}B^TS(t)\mathbf{x}(t) \tag{3}$$

where

$$\dot{S}(t) = -S(t)A - A^TS(t) + S(t)BR^{-1}B^TS(t) - Q; \quad S(t_\mathrm{f}) = S_\mathrm{f} \tag{4}$$

Here, $S(t)$ denotes the time-varying Riccati matrix. It is assumed that $(A, B)$ is stabilizable, and $(A, Q^{1/2})$ is detectable.

## WEIGHT OPTIMIZATION USING CLOSED-FORM ALGEBRAIC SOLUTIONS

Classically Eq. (4) is numerically integrated backward in time from $t = t_\mathrm{f}$ to $t = t_0$ to compute the controlled system response forward in time using the optimal control defined by introducing Eq. (3) into Eq. (2). After this process, the controlled response is optimized by assuming that the symmetric free elements of the weight matrices are available for optimization. A Taylor expansion model is developed for the augmented state consisting of the controlled state and controls. Closed-form algebraic models are developed for the Taylor expansion partial derivatives that are required for the free weight matrix parameters. These equations are generated by computing the partial derivatives for Eqs. (2) and (4), and the state and Riccati matrix sensitivity partial derivatives are derived as

$$\dot{\mathbf{x}}_{,p}(t) = A\mathbf{x}_{,p}(t) + B\mathbf{u}_{,p}(t); \quad \mathbf{x}_{,p}(t_0) = 0 \tag{5}$$

$$\begin{aligned}\dot{S}_{,p}(t) = &-S_{,p}(t)A - A^T S_{,p}(t) + S_{,p}(t)BR^{-1}B^T S(t) + S(t)BR^{-1}B^T S_{,p}(t) \\ &- S(t)BR^{-1}R_{,p}R^{-1}B^T S(t) - Q_{,p}; \quad S_{,p}(t_0) = 0\end{aligned} \tag{6}$$

where $(\cdot)_{,p}$ represents the partial derivative with respect to each symmetric element of $Q$, $R$, and $S_\text{f}$. The number of symmetric elements is $M = 2n_\text{sym} + m_\text{sym}$, where is $n_\text{sym} = n(n+1)/2$ and $m_\text{sym} = m(m+1)/2$. The sensitivity equations are simultaneously integrated with the state and Riccati matrix differential equations. As a result, $M + 1$ state and state partial equations and $M + 1$ Riccati Matrix and Riccati Matrix partial equations are numerically integrated. Though technically straightforward, as the state dimension increases, the computational impact on high-order systems is very significant. This paper presents a novel computational strategy that eliminates the requirement for invoking the computationally intense numerical integration step for the state, Riccati matrix, and associated partial derivatives. Closed-form algebraic equations are developed for all key equations and their partial derivatives.

## Closed-Form Solution for the Time-Varying Riccati Matrix, the State Trajectory, and the Performance Index

### *Closed-Form Solution for the Time-Varying Riccati Matrix*

Numerical integrations for Eqs. (4) and (6) are eliminated by introducing the following closed-form solution for the differential matrix Riccati equation, where the solution consists of steady-state and time-varying parts[14,15]

$$S(t) = S_\text{ss} + Z^{-1}(t); \quad S(t_\text{f}) = S_\text{f} \tag{7}$$

where the steady-state solution for $S_{ss}$ satisfies the algebraic matrix Riccati equation:

$$0 = -S_\text{ss}A - A^T S_\text{ss} + S_\text{ss}BR^{-1}B^T S_\text{ss} - Q \tag{8}$$

By substituting Eq. (7) into Eq. (4), the matrix differential Lyapunov equation for $Z(t) \in \Re^{n\times n}$ is derived as

$$\dot{Z}(t) = \bar{A}Z(t) + Z(t)\bar{A}^T - BR^{-1}B^T; \quad Z(t_\text{f}) = (S_\text{f} - S_\text{ss})^{-1} = Z_\text{b} \tag{9}$$

where $\bar{A} = A - BR^{-1}B^T S_\text{ss}$ is the system stability matrix.[16] The closed-form solution for $Z(t)$ consists of a steady-state and time-varying parts, leading to

$$Z(t) = Z_\text{ss} + e^{\bar{A}(t-t_\text{f})}[(S_\text{f} - S_\text{ss})^{-1} - Z_\text{ss}]\, e^{\bar{A}^T(t-t_\text{f})} \tag{10}$$

where $e^{(\cdot)}$ is the $\Re^{n\times n}$ exponential matrix,[17] $S_\text{f} \neq S_\text{ss}$, and the steady-state solution for $Z_\text{ss}$ satisfies the algebraic matrix Lyapunov equation:[18]

$$\bar{A}Z_\text{ss} + Z_\text{ss}\bar{A}^T - BR^{-1}B^T = 0 \tag{11}$$

### *Closed-Form Solution for the State Trajectory*

Introducing Eq. (3) into Eq. (2), the governing differential equation for the controlled state follows as

$$\dot{\mathbf{x}}(t) = [A - BR^{-1}B^T S(t)]\mathbf{x}(t); \quad \mathbf{x}_0 = \mathbf{x}(t_0) \tag{12}$$

where $S(t)$ is given by Eq. (4). Turner and Chun[19,20] have shown that the closed-form solution for Eq. (12) is given by the mapping equation:

$$\mathbf{x}(t) = \phi(t, t_0)\mathbf{x}_0 \tag{13}$$

where $\phi(t, t_0)$ denotes the state transition matrix, which has the explicit form:

$$\phi(t, t_0) = Z(t)e^{-\bar{A}^T(t-t_0)}Z^{-1}(t_0) \tag{14}$$

where $Z(t)$ is defined by Eq. (10), and $\phi(t, t_0)$ satisfies the standard group properties: $\phi(t_2, t_0) = \phi(t_2, t_1)\phi(t_1, t_0)$ and $\phi^{-1}(t_1, t_0) = \phi(t_0, t_1)$.

*Closed-Form Solution for the Performance Index*

The performance index cost can be expressed as:

$$J = \frac{1}{2}\mathbf{x}(t)^T S(t)\mathbf{x}(t) \tag{15}$$

where $S(t)$ is defined by Eq. (4), $\mathbf{x}(t)$ is defined by Eq. (13), and $Z(t)$ is defined by Eq. (10). It is important to observe that no numerical integration is required because all equations are in a closed-form. Figure 1 displays the difference between the conventional and proposed approaches for solving the optimal feedback control problem. As mentioned earlier, the proposed approach eliminates the numerical integration process that is computationally intensive by introducing closed-form solutions for each variable.

| Conventional approach | Proposed approach |
|---|---|
| Given the optimal feedback control problem<br>$J = \frac{1}{2}\mathbf{x}^T(t_f)S_f\mathbf{x}(t_f) + \frac{1}{2}\int_{t_0}^{t_f}(\mathbf{x}^T(t)Q\mathbf{x}(t) + \mathbf{u}^T(t)R\mathbf{u}(t))\,dt$<br>Subject to $\dot{\mathbf{x}}(t) = A\mathbf{x}(t) + B\mathbf{u}(t)$; $\mathbf{x}(t_0) = \mathbf{x}_0, \mathbf{x}(t_f) = \mathbf{0}$. | |
| Performing numerical integration backward in time<br>$\dot{S}(t) = -S(t)A - A^T S(t) + S(t)BR^{-1}B^T S(t) - Q$;<br>$S(t_f) = S_f$ | $S(t) = S_{ss} + Z^{-1}(t)$; $S(t_f) = S_f$<br>$0 = -S_{ss}A - A^T S_{ss} + S_{ss}BR^{-1}B^T S_{ss} - Q$<br>Solving Riccati equation → $S_{ss}$<br>$\bar{A}Z_{ss} + Z_{ss}\bar{A}^T - BR^{-1}B^T = 0$, where $\bar{A} = A - BR^{-1}B^T S_{ss}$<br>Solving Lyapunov equation → $Z_{ss}$<br>$Z(t) = Z_{ss} + e^{\bar{A}(t-t_f)}[(S_f - S_{ss})^{-1} - Z_{ss}]\,e^{\bar{A}^T(t-t_f)}$ |
| $\mathbf{x}(t) = e^{A-BR^{-1}B^T S(t)(t-t_0)}\mathbf{x}_0$<br>$\mathbf{u}(t) = -R^{-1}B^T S(t)\mathbf{x}(t)$<br>$J = \frac{1}{2}\mathbf{x}(t)^T S(t)\mathbf{x}(t)$ | $\mathbf{x}(t) = Z(t)e^{-\bar{A}^T(t-t_0)}Z^{-1}(t_0)\mathbf{x}_0$<br>$\mathbf{u}(t) = -R^{-1}B^T S(t)\mathbf{x}(t)$<br>$J = \frac{1}{2}\mathbf{x}(t)^T S(t)\mathbf{x}(t)$ |

**Figure 1. Comparison of Solving the Optimal Feedback Control Problem**

## Optimization of the Weight Matrices

The unknown parameters for the optimization process are the symmetric elements of $Q$, $R$, and $S_f$. To optimize the components of the weight matrices, analytic gradients are computed for the optimization steps. The optimization process finds the weight matrices to minimize the terminal values of the state and control to be zero.

The goal of the weight matrices optimization is to develop a state-space model consisting of the state and control, which is called an augmented state as follows:

$$\mathbf{y}(t) = [\mathbf{x}^T(t), \mathbf{u}^T(t)]^T \tag{16}$$

which is evaluated at the final time for the optimization process. Computationally the optimization seeks to force the terminal values of $\mathbf{y}(t)$ to zero. Equation (16) is Taylor expanded at $t = t_\mathrm{f}$ as a function of the weight matrices' symmetric elements, which leads to

$$\mathbf{y}_\mathrm{f} = \mathbf{y}_\mathrm{f}(Q, R, S_\mathrm{f}) + \sum_k^{n_\mathrm{sym}} \frac{\partial \mathbf{y}}{\partial q_k} \mathrm{d}q_k + \sum_l^{m_\mathrm{sym}} \frac{\partial \mathbf{y}}{\partial r_l} \mathrm{d}r_l + \sum_k^{n_\mathrm{sym}} \frac{\partial \mathbf{y}}{\partial s_{\mathrm{f}_k}} \mathrm{d}s_{\mathrm{f}_k} \tag{17}$$

Assuming that the left-hand side of Eq. (17) is zero (as mentioned earlier for derivation simplification purposes) and collecting the partials into a global Jacobian matrix $\Delta J$, leads to

$$0 = \mathbf{y}_\mathrm{f}(Q, R, S_\mathrm{f}) + \Delta J \mathrm{d}\mathbf{p} \tag{18}$$

where

$$\Delta J = \left[\frac{\partial \mathbf{y}}{\partial q_1}, \cdots, \frac{\partial \mathbf{y}}{\partial q_{n_\mathrm{sym}}}, \quad \frac{\partial \mathbf{y}}{\partial r_1}, \cdots, \frac{\partial \mathbf{y}}{\partial r_{m_\mathrm{sym}}}, \quad \frac{\partial \mathbf{y}}{\partial s_1}, \cdots, \frac{\partial \mathbf{y}}{\partial s_{n_{sym}}}\right] \tag{19}$$

$$\mathrm{d}\mathbf{p} = \left[\mathrm{d}q_1, \cdots, \mathrm{d}q_{n_\mathrm{sym}}, \quad \mathrm{d}r_1, \cdots, \mathrm{d}r_{m_\mathrm{sym}}, \quad \mathrm{d}s_1, \cdots, \mathrm{d}s_{n_{sym}}\right]^T \tag{20}$$

The sensitivity partial derivatives used in the optimization process are discussed in the following subsection.

*Minimum Norm Correction Solution*

Since the number of unknowns is larger than the number of end conditions, the solution for the parameter correction vector $\mathrm{d}\mathbf{p}$ is defined by minimizing the following performance index:

$$\Phi = \frac{\mathrm{d}\mathbf{p}^T \mathrm{d}\mathbf{p}}{2} + \boldsymbol{\lambda}^T [\mathbf{y}_\mathrm{f} + \Delta J \mathrm{d}\mathbf{p}] \tag{21}$$

where $\boldsymbol{\lambda}$ denotes the Lagrange multiplier, and the necessary condition for the optimization follows as

$$\begin{gathered} \Phi_{,\mathrm{d}\mathbf{p}} = \mathrm{d}\mathbf{p} + \Delta J^T \boldsymbol{\lambda} = \mathbf{0} \\ \Phi_{,\boldsymbol{\lambda}} = \mathbf{y}_\mathrm{f} + \Delta J \mathrm{d}\mathbf{p} = \mathbf{0} \end{gathered} \tag{22}$$

The solution for $\boldsymbol{\lambda}$ is given by

$$\boldsymbol{\lambda} = (\Delta J \Delta J^T)^{-1} \mathbf{y}_\mathrm{f} \tag{23}$$

and the minimum norm optimization solution for the parameter correction vector is obtained as

$$\mathrm{d}\mathbf{p} = -\Delta J^T (\Delta J \Delta J^T)^{-1} \mathbf{y}_\mathrm{f} \tag{24}$$

Therefore, the elements of the weight matrices are updated by the following equation:

$$\mathbf{p}_\mathrm{update} = \mathbf{p}_\mathrm{previous} + \mathrm{d}\mathbf{p} \tag{25}$$

where the vector $\mathbf{p}$ is defined as

$$\mathbf{p} = \left[q_1, \ldots, q_{n_{sym}}, r_1, \ldots, r_{m_{sym}}, s_1, \ldots, s_{n_{sym}}\right]^T \tag{26}$$

The entire procedure of the optimization process is described in Figure 2. At the beginning of the optimization process, initial weight matrices are selected, and $\mathbf{p}$ that consists of the symmetric elements of weight matrices is defined. After computing the Jacobian matrix and the augmented states at the final time, $\mathrm{d}\mathbf{p}$ is calculated using Eq. (24). Then, the new symmetric elements of weight matrices are updated using Eq. (25). Here, the weight matrices should satisfy the definiteness conditions as mentioned in the problem formulation. Hence, after checking the definiteness of each weight matrix, the weight matrices that only satisfy the definiteness conditions are updated and conveyed to the next iteration. If the solution of the optimal feedback control problem using the updated weight matrices does not meet the requirement defined by users, one performs the next iteration step to update the weight matrices.

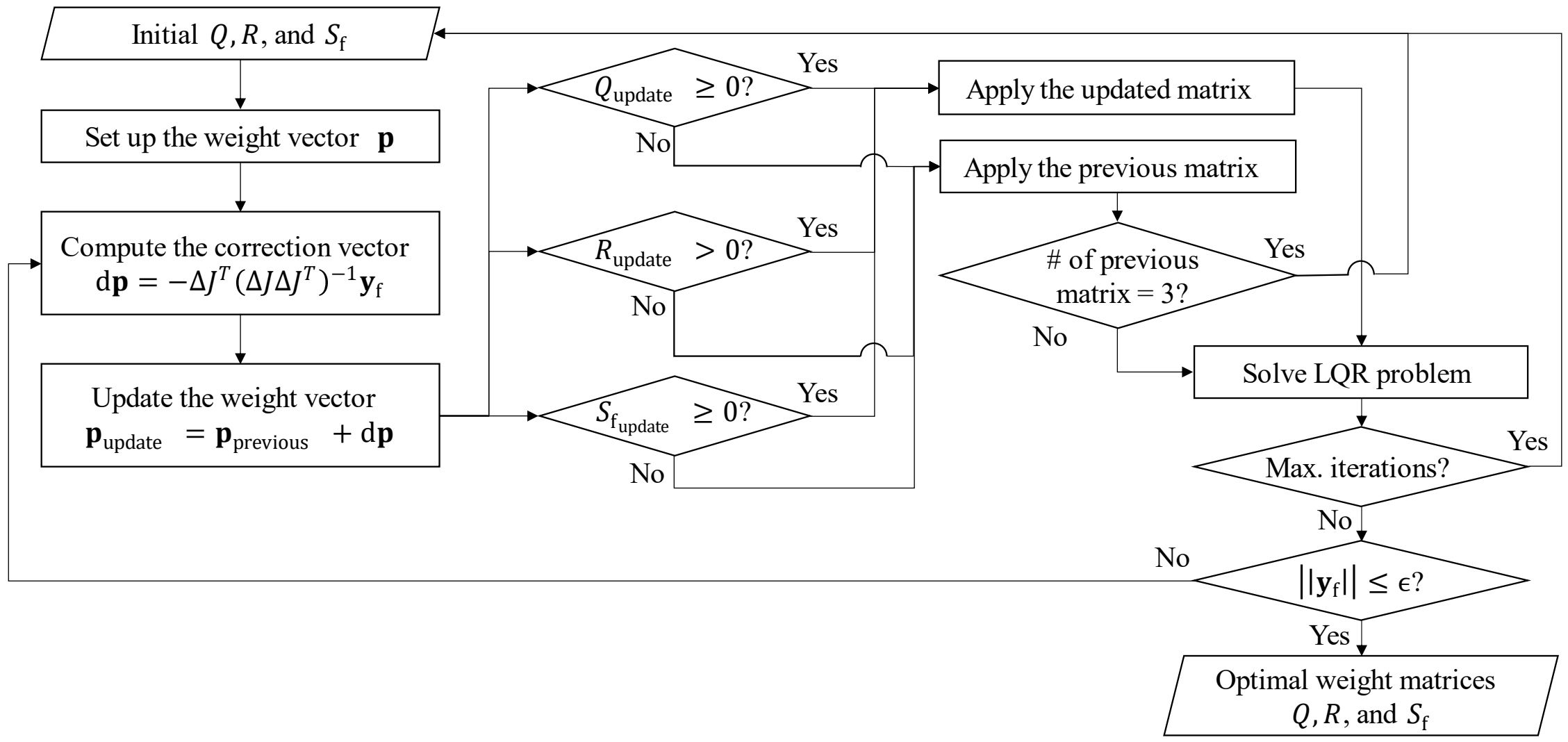

**Figure 2. Flowchart of the Optimization Process for the Weight Matrices**

### Sensitivity Partial Derivatives for the State, Control, and Riccati Matrix

To form the Jacobian matrix in Eq. (19), the sensitivity partial derivatives for the state, control, and Riccati matrix are required. If Eqs. (5) and (6) are directly used and numerically integrated to form the Jacobian matrix, then it requires a high computational burden when increasing the dimension of the states. For this reason, this work presents partial derivatives with closed-form solutions for all key equations.

*State Partials*

The state trajectory partial derivatives are generated by modeling the state as Eq. (13), leading to the state partial derivatives given by

$$\mathbf{x}_{,p}(t) = \phi_{,p}(t, t_0)\mathbf{x}_0 \tag{27}$$

where the state transition matrix partial derivatives are given by

$$\begin{aligned}\phi_{,p}(t, t_0) = Z_{,p}(t)e^{-\bar{A}^T(t-t_0)}Z^{-1}(t_0) + Z(t)\left[e^{-\bar{A}^T(t-t_0)}\right]_{,p} Z^{-1}(t_0) \\ - Z(t)e^{-\bar{A}^T(t-t_0)}Z^{-1}(t_0)Z_{,p}(t)Z^{-1}(t_0)\end{aligned} \tag{28}$$

Here, $[\cdot]_{,p}$ denotes the partial derivative of the matrix exponential. The main point is that all of the equations in Eq. (28) are purely algebraic.

*Control Partials*

The partial derivatives of the feedback control of Eq. (3) are expressed as

$$\begin{aligned}\mathbf{u}_{,p}(t) = R^{-1}R_{,p}R^{-1}B^TS(t)\phi(t, t_0)\mathbf{x}_0 - R^{-1}B^TS_{,p}(t)\phi(t, t_0)\mathbf{x}_0 \\ - R^{-1}B^TS(t)\phi_{,p}(t, t_0)\mathbf{x}_0\end{aligned} \tag{29}$$

and Eqs. (27) and (29) are purely algebraic in nature. That is, no numerical integration is required.

*Riccati Matrix partial derivatives*

From Eq. (7), the partial derivatives of the Riccati matrix are given by

$$S_{,p}(t) = S_{\mathrm{ss},p} - Z^{-1}(t) Z_{,p}(t) Z^{-1}(t) \tag{30}$$

The partial derivatives of the time-varying part of the Riccati matrix are expressed as

$$\begin{aligned} Z_{,p}(t) = Z_{\mathrm{ss},p} + \left[e^{\bar{A}(t-t_\mathrm{f})}\right]_{,p} Z_\mathrm{b} e^{\bar{A}^T(t-t_\mathrm{f})} + e^{\bar{A}(t-t_\mathrm{f})} Z_{\mathrm{b},p}\, e^{\bar{A}^T(t-t_\mathrm{f})} \\ + e^{\bar{A}(t-t_\mathrm{f})} Z_\mathrm{b} \left[e^{\bar{A}^T(t-t_\mathrm{f})}\right]_{,p} \end{aligned} \tag{31}$$

where $Z_\mathrm{b}$ is the boundary condition for $Z(t)$, and the partial derivatives of $Z_\mathrm{b}$ in Eq. (9) are given by

$$Z_{\mathrm{b},p} = (S_\mathrm{f} - S_\mathrm{ss})^{-1} S_{\mathrm{ss},p} (S_\mathrm{f} - S_\mathrm{ss})^{-1} - Z_{\mathrm{ss},p} \tag{32}$$

In addition, the partial derivatives for the steady-state Riccati matrix follow as the matrix Lyapunov equation:

$$\bar{A}^T S_{\mathrm{ss},p} + S_{\mathrm{ss},p} \bar{A} = -S_\mathrm{ss} B R^{-1} R_{,p} R^{-1} B^T S_\mathrm{ss} - Q_{,p} \tag{33}$$

where $\bar{A} = A - BR^{-1}B^T S_\mathrm{ss}$. The partial derivatives for the steady-state Lyapunov matrix follow as

$$\bar{A} Z_{\mathrm{ss},p} + Z_{\mathrm{ss},p} \bar{A}^T = -\bar{A}_{,p} Z_\mathrm{ss} - Z_\mathrm{ss} \bar{A}^T{}_{,p} - BR^{-1} R_{,p} R^{-1} B^T \tag{34}$$

where the closed-loop system dynamics matrix partial derivatives are given by

$$\bar{A}_{,p} = BR^{-1} R_{,p} R^{-1} B^T S_\mathrm{ss} - BR^{-1} B^T S_{\mathrm{ss},p} \tag{35}$$

## SIMULATION STUDY

To validate the performance of the weight optimization process, an example problem is considered for the second-order differential equation. The optimal feedback control problem is defined in Eqs. (1) and (2), and the system dynamics matrix and control influence matrices are defined as

$$\mathbf{x}(t) \in \mathbb{R}^2, u(t) \in \mathbb{R}, A = \begin{bmatrix} 0 & 1 \\ -k/m & -c/m \end{bmatrix}, B = \begin{bmatrix} 0 \\ m \end{bmatrix} \tag{36}$$

In addition, the parameters used in the simulations are listed in Table 1.

**Table 1. Simulation Parameters**

| Variables | Values |
|---|---|
| Mass (kg) | $m = 1$ |
| Spring coefficient (N/m) | $k = 0.64$ |
| Damping coefficient (Ns/m) | $c = 0.16$ |
| Initial condition $\mathbf{x}(t_0)$ | $[10 \quad 10]^T$ |
| Final time $t_\mathrm{f}$ (sec) | 10 |
| Time interval $\mathrm{d}t$ (sec) | 0.01 |

To optimally determine the weight matrices, the weight optimization process proposed is validated in this section. The goal of the optimization process is to find the optimal weight matrices that make the states and the control input at the final time close to zero, minimizing the corrections for the weight matrices. The requirement is set to be the norm of the augmented state at the final time that is less than $10^{-5}$.

**Table 2. Optimization Results**

| | Initial | Optimized |
|---|---|---|
| Norm of the augmented state $\|\|\mathbf{y}_\mathrm{f}\|\|$ | $2.68 \times 10^{-2}$ | $3.82 \times 10^{-6}$ |
| States at the final time $\mathbf{x}(t_\mathrm{f})$ | $[1.20, -1.70]^T \times 10^{-2}$ | $[3.28, -1.94]^T \times 10^{-6}$ |
| Control at the final time $\mathbf{u}(t_\mathrm{f})$ | $1.70 \times 10^{-2}$ | $-2.71 \times 10^{-7}$ |
| State weight matrix $Q$ | $I_{2\times 2}$ | $\begin{bmatrix} 1.63 & 0.08 \\ 0.08 & 1.02 \end{bmatrix}$ |
| Control weight matrix $R$ | 1 | 0.19 |
| Final state weight matrix $S_\mathrm{f}$ | $I_{2\times 2}$ | $\begin{bmatrix} 1.08 & 0.31 \\ 0.31 & 0.50 \end{bmatrix}$ |
| Performance index $J$ | $2.01 \times 10^{2}$ | $1.58 \times 10^{2}$ |

Table 2 shows the optimization results when the initial weight matrices are assumed as identity matrices. The norm of the augmented state at the final time using initial weight matrices does not satisfy the requirement, but this value reaches $3.82 \times 10^{-6}$, which satisfies the requirement, using the proposed approach. In Figure 3, it is shown that the norm of the augmented state is linearly decreased in the log scale over the iterations. After the optimization process, the performance index is also decreased compared to one using the initial weight matrices. Furthermore, as shown in Figure 4, the proposed approach optimizes the off-diagonal terms, as well as the diagonal terms, albeit the initial weight matrices are set as identity matrices only. Also, all the weight matrices obtained are positive definite. By optimizing all components of the weight matrices, it provides more flexibility for the optimization. Figures 5 and 6 show the states and optimal control trajectories obtained from the optimization process for all iterations, and the response of the states and control becomes faster over the iterations to meet the requirement.

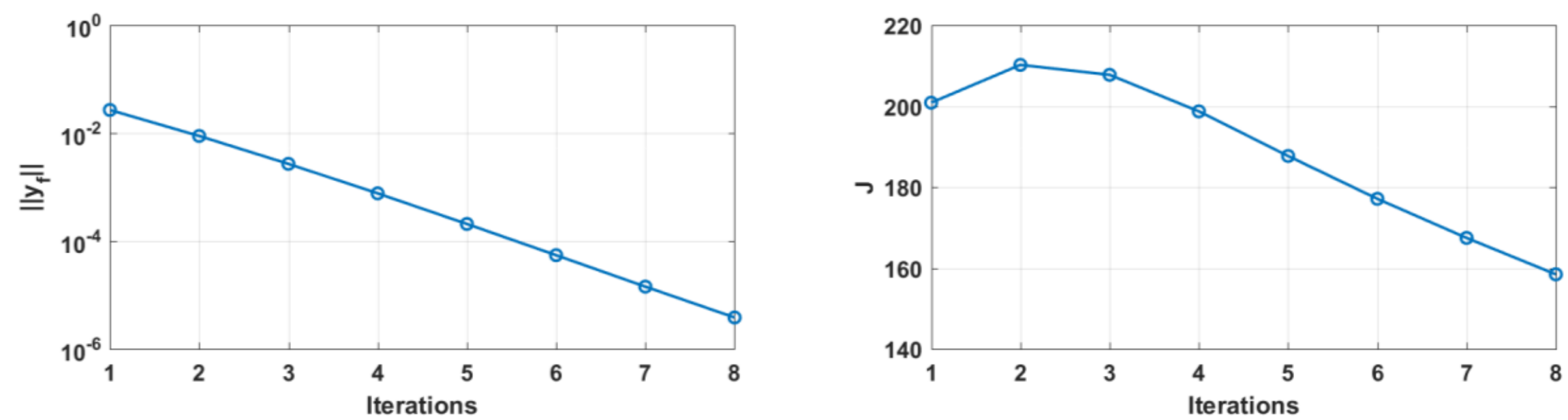

**Figure 3. History of the Norm of the Augmented State (Left) and the Performance Index (Right)**

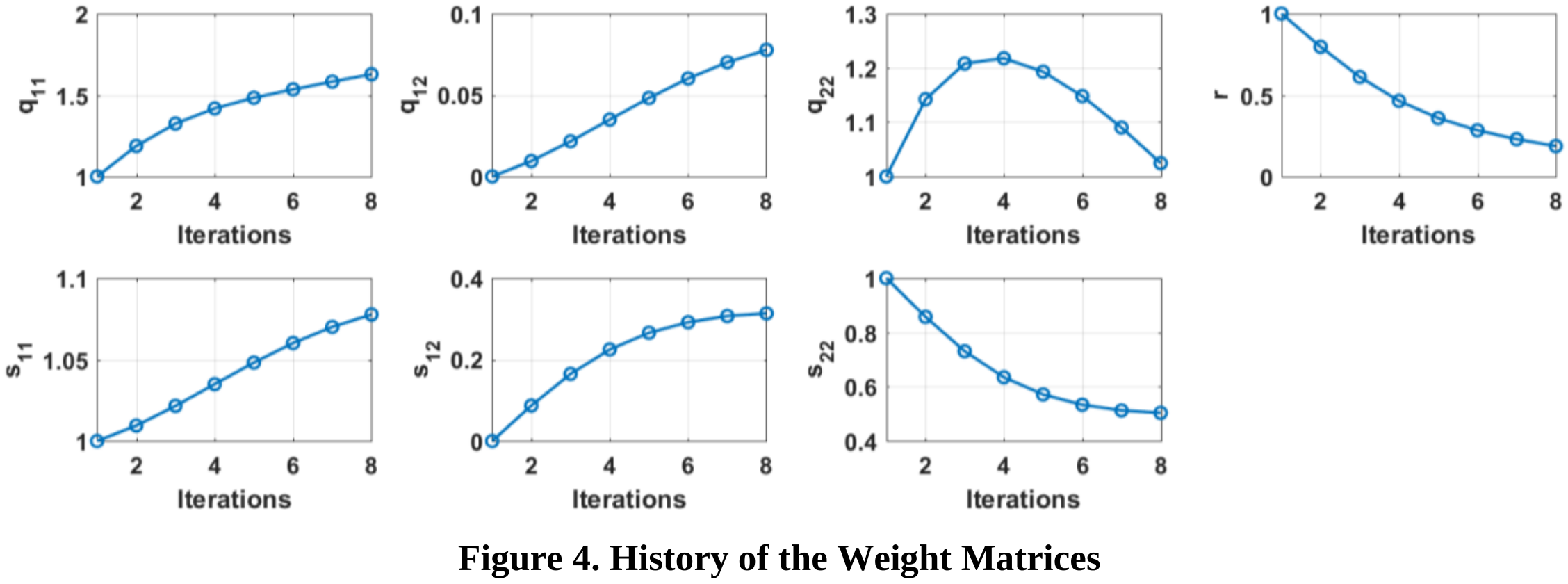

**Figure 4. History of the Weight Matrices**

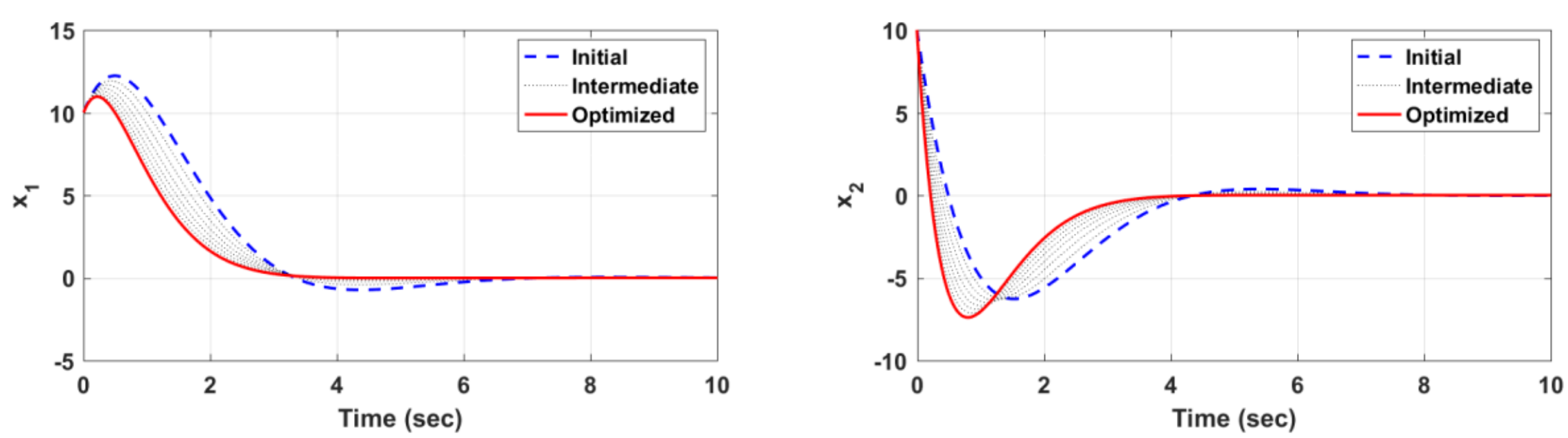

**Figure 5. State Trajectories**

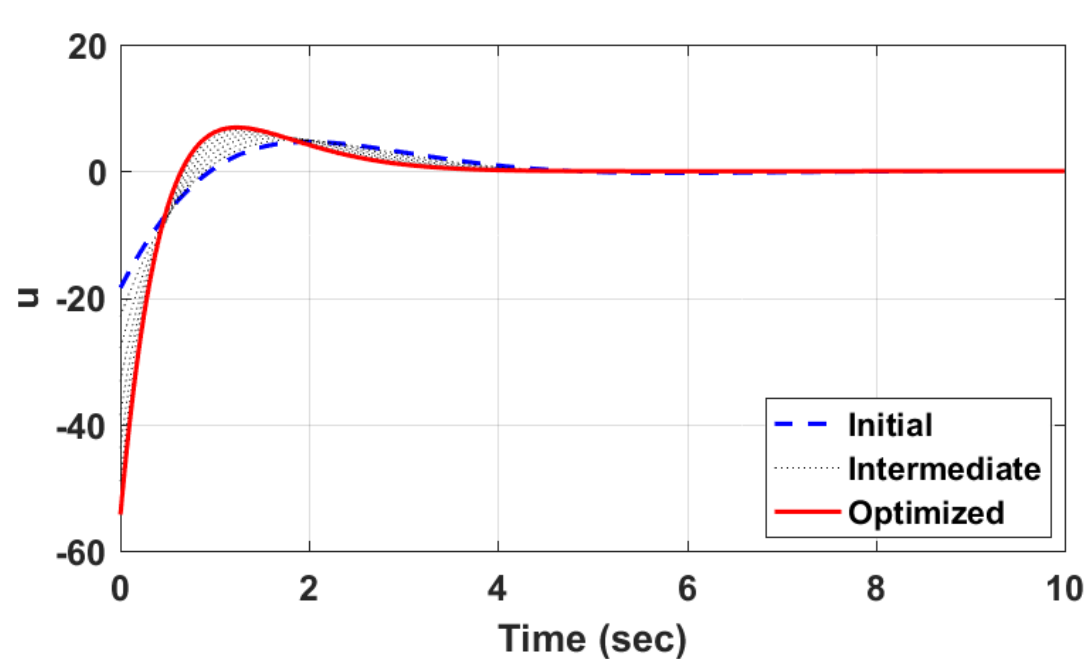

**Figure 6. Optimal Control Trajectories**

## CONCLUSION

This work proposes an optimization approach for the weight matrices of the feedback control problem with a quadratic performance index. The proposed optimization process is defined by Taylor's expansion of the state and control at the final time to force the desired final time boundary conditions for the system response. The optimization process contains multiple numerical integrations for the states, differential Riccati equations, and their partial derivatives. The volume of the calculations increases as the dimension of the states increases. To reduce the computational burden of the optimization process, this work utilizes the algebraic equations that exploit the closed-form solutions for the time-varying states, state transition matrix, and matrix exponential as well as the time-varying and algebraic solutions for Lyapunov and Riccati matrix equations. The proposed weight matrices optimization approach utilizing the closed-form solutions is validated by numerical simulations for the dynamic system. Consequently, the proposed approach provides the optimal weight matrices that satisfy the required norm of the augmented states that contains the states and control at the final time. The proposed computational procedures are expected to the broadly useful for control theory applications in science and engineering fields.